\documentclass[11pt]{amsart}
\usepackage{amsmath}
\usepackage{amssymb}
\usepackage{mathrsfs}
\usepackage{tikz}
\newcommand{\la}{\langle}
\newcommand{\ra}{\rangle}
\newtheorem{teo}{Theorem}
\newtheorem{pro}{Proposition}
\newtheorem{cor}{Corollary}
\newtheorem{lem}{Lemma}
\newtheorem*{nota}{Remark}
\newtheorem*{note}{Remarks}

\title{Reflected Brownian motion in Weyl chambers}

\begin{document}
\maketitle
\centerline{N. DEMNI \footnote{supported by CMCU 07G1501 , e-mail: demni@math.uni-bielefeld.de}}

\begin{abstract}
We supply two different descriptions of the pushing process driving the reflected Brownian motion in Weyl chambers, when the latter domains are simplexes. The first one shows that a simple root lies in one and only one orbit if and only if the pushing process in the direction of that simple root increases as the sum of all the Brownian local times in the directions of the orbit's positive elements. The last one shows that the pushing process may be written as the sum of an inward normal vector at the chamber's boundary and an inward normal vector at the origin, yielding a kind of a multivoque stochastic differential equation for the reflected process. We finally give a particles system interpretation of the reflected process and construct a multidimensional  skew Brownian motion.  
\end{abstract} 

\section{overview}
By `the reflected Brownian motion' (shorthand RBM), it is often meant the absolute value of a real Brownian motion. This process has gained much fame since it is not an It\^o's semimartingale and due to the celebrated L\'evy's representation in relation with Skorohod's problem (\cite{Rev}). More precisely, if $B$ is a standard real Brownian motion, then there exist a standard real Brownian motion $\beta$  and an increasing process  $L^0(B)$ such that 
 \begin{equation}\label{E1}
|B_t| = \beta_t + L_t^0(B). 
\end{equation}
The process $L^0$ is known as the local time at $0$ of $B$ since the following representation holds (Ch.VI in \cite{Rev})
\begin{equation*}
L_t^0(B) = \lim_{\epsilon \rightarrow 0} \frac{1}{2\epsilon} \int_0^t {\bf 1}_{\{|B_s| \leq \epsilon\}} ds = \frac{1}{2}L_t^0(|B|),
\end{equation*}
and it gives the explicit Doob-Meyer's decomposition of the submartingale $|B|$. It also shows that $L^0(B)$ is adapted with respect to the natural filtration $\mathscr{F}^{|B|}$ of $|B|$. Besides, the support of $dL^0(B)$ is contained in the set $\{t, B_t = 0\}$ and the identity in law of trajectories due to Paul L\'evy holds (\cite{Rev})
\begin{equation*}
(|B|, L^0(B)) \overset{d}{=} (S - B, S),\quad S_t := \sup_{0 \leq s \leq t}B_s.
 \end{equation*}
Motivated by applications to queueing theory and constructions of stable subordinators, various probabilists were interested in multidimensional analogs of the RBM. Such processes behave as multidimensional Brownian motions in the interior of their state spaces and perform reflections at the boundaries. Examples of state spaces already appeared in literature include  wedges (\cite{Var}), two dimensional cones 
\begin{equation*}
\{z = re^{i\theta} \in \mathbb{C},\, |\theta| \leq \alpha, \, 0 < \alpha \leq \pi/2\}
\end{equation*} and more generally simplexes 
\begin{equation*}
 \{x \in \mathbb{R}^n, x = \sum_{i=1}^n x_i v_i, x _i \geq 0\}
\end{equation*}
with respect to some basis $(v_i)_i$ so that a side, say $\Delta_i, 1 \leq i \leq n$ is specified by $x_i = 0$ (\cite{Legall}), domains delimited by a set of hyperplanes (\cite{Sznit}), convex polyhedral domains: 
\begin{equation*}
\{x \in \mathbb{R}^n,\, \la n_i, x \ra \geq b_i, \, 1 \leq i \leq d\}  
\end{equation*}
where $(n_i)_{1 \leq i \leq d}$ are $n$-dimensional vectors and $(b_i)_{1 \leq i \leq d}$ is a $d$-dimensional vector  (\cite{Dai}) and in particular orthants (\cite{Tay}). Most of the papers cited above were mainly intended to decide whether the RBM is a semimartingale or not and whether the uniqueness in law holds in the affirmative case. \\ 
Recently, O. Chybiryakov defined and studied a RBM valued in a class of cones that arises in Lie  algebra theory and more precisely associated with root systems:  Weyl chambers (\cite{Chy},\cite{Hum}). From their very definition (see below), these domains are convex simple polyhedrons for which $(n_i)_i$ are the so-called simple roots and $b_i = 0$ for all $i$, are simplexes for most of the root systems, in particular are orthants for the so-called orthogonal root system and are wedges of angles $\pi/m, m \in \{2,3\dots \}$ for dihedral systems (\cite{Hum}). When it is not a simplex, the Weyl chamber may be a half plane (root systems of type $A_2$), a setting that was of considerable interest to various mathematicians (see p.406 in \cite{Var} for references). Assuming that Weyl chambers are simplexes, the author actually derived a Tanaka's type SDE and a representation of its drift, the pushing process, as an additive functional.  With regard to the wealthy and regular geometrical structure of Weyl chambers and motivated by applications to queueing theory, we thought it is natural to answer the following questions:  
\begin{itemize}
\item Given a simple root, is it possible to write the coordinate of the pushing process in this direction as a sum of local times of some real semimartingales?
\item Is there any relation to multivoque SDEs, a tool that showed to be powerful in stochastic analysis with reflections (\cite{Cepa})?
\end{itemize}
The answer to the first question is positive if and only if no other simple root lies in the orbit of the given simple root, while there is no perfect relation to multivoque stochastic analysis and this fact is not surprising since the reflections are not normal to the boundary. Nevertheless, we show that the additive functional decomposes as the sum of an inward normal vector at the boundary of the Weyl chamber and of an inward normal vector at the vertex of the simplex. As a by product, the RBM satisfies a multivoque-type SDE with a repelling to infinity which agrees with the transience of the reflected process. The paper is closed by further developments: we express the pushing process by means of local times of the RBM and we give a construction of a multidimensional the skew Brwonian motion (\cite{Lejay}). For sake of completeness, we collect below some facts on root systems and fix some notations we use later on (Ch.VI in \cite{Dunkl}, Ch.I, II in \cite{Hum}). Once we do, we state Chybiryakov's non published result and we review its proof. This is two-fold: we correct some minor errata and we discuss the case of the root system of type $A$ for which $\overline{C}$ is not a simplex of $V$.

\section{RBM in Weyl chambers: Reminder}

\subsection{Facts on root systems}  Let  $(V, \la,\ra)$ be a $n$-dimensional Euclidean space and define the reflection 
\begin{equation*}
\sigma_{\alpha}(v) = v - 2\frac{\la \alpha,v\ra}{\la \alpha, \alpha \ra} \alpha, \quad \alpha \in V \setminus \{0\} 
\end{equation*}
with respect to the hyperplane $H_{\alpha} = \alpha^{\perp} := \{x, \la \alpha, x \ra = 0\}$. Then, a \emph{root system} $R$ is a collection $\{\alpha \in V\setminus \{0\}\}$ such that $\sigma_{\alpha}(R) = R$ for all $\alpha \in R$. It is said to be reduced if  $\mathbb{R}\alpha \cap R = \{\pm \alpha\}$. A trivial example is provided by the {\it orthogonal} root system 
$R = \{\pm e_i, 1 \leq i \leq n\}$, where $(e_i)_i$ is the canonical basis of $V$. It can be shown that there exists a basis $S$ of $\textrm{span}(R) \subset V$, called simple system, such that a root $\alpha \in R$ is either a positive or a negative linear combination of elements of $S$. Thus $R$ splits into a positive and a negative systems $R_+, R_-$ respectively. With the help of elements of $S$, known as simple roots, one defines the so-called \emph{closed Weyl chamber}  
\begin{equation*}
\overline{C} := \{x \in V, \, \la \alpha, x \ra \geq 0 \, \textrm{for\, all}\, \alpha \in S\}.
\end{equation*}  
This is a cone, not necessarily a simplex, and specializes to the positive orthant of $\mathbb{R}^n$ when $R$ is the orthogonal root system. Moreover, $\overline{C}$ enjoys the following important property: let $W$ be the group spanned by $\{\sigma_{\alpha}, \alpha \in R\}$: the reflections group (it is a subgroup of the orthogonal group $O(n)$). Then, $\overline{C}$ is a {\it fundamental domain} with respect to the action of $W$ on $V$ ($x \mapsto wx$), in the sense that for any $x \in V$, there exist $w \in W$ and a unique $y \in \overline{C}$ such that $x = wy$.  Thus, the space $V$ decomposes as 
\begin{equation*}
V = \cup_{w \in W} w \overline{C} = \left(\cup_{w \in W} w C\right) \cup \left(\cup_{\alpha \in S}H_{\alpha}\right)
\end{equation*}
and the projection $\pi: V \mapsto V^W \mapsto \overline{C}$ is well defined (note that $\pi = w^{\star}$ if $x \in w\overline{C}$ for some $w \in W$), where $V^W$ is the orbits space.
When $S = \{s_1,\dots, s_n\}$ is a basis of $V$ then $\overline{C}$ is a simplex of vertex $0_V$ since
\begin{equation*}
\overline{C} = \{x \in V,\, x = \sum_{i=1}^n \la s_i,x \ra \xi_i, \la s_i, x \ra \geq 0,\,1 \leq i \leq n\},  
\end{equation*}
where  $(\xi_i)_i$ is the dual basis of $S$ defined by: 
\begin{equation*}
\la \xi_i,s_j \ra = \delta_{ij}.
\end{equation*} 
In the last decade, an active probabilistic research related to root systems has emerged after the introduction of the so-called Dunkl processes and their $W$-invariant parts (\cite{Chy1}). The latter processes are diffusions valued in $\overline{C}$ and specialize to the RBM in Weyl chambers defined as $\pi(\Theta)$, where $\Theta$ is a $V$-valued Brownian motion.  

\subsection{On Chybiryakov's result and review of its proof}
\begin{teo}\label{T1}
Assume that $S = \{s_1, \dots, s_n\}$ is a basis of $V$. Then there exist a $V$-valued Brownian motion $B$ and a $\overline{C}$-valued continuous process denoted $L^0(\Theta)$  such that
\begin{equation}\label{SDE}
\pi(\Theta_t) = \pi(\Theta_0) + B_t + L_t^0(\Theta),
\end{equation}  
and for any $1 \leq i \leq n$, $\la s_i, L^0(\Theta) \ra$  is an increasing process satisfying 
\begin{equation}\label{Local}
\int_0^t {\bf 1}_{\{\la s_i , \pi(\Theta_s)\ra \neq 0\}} \la s_i, dL_s^0(\Theta)\ra = 0.  
\end{equation}
More precisely, if $\Lambda$ is the matrix whose rows are the simple roots, then 
\begin{equation*}
\Lambda L^0(\Theta) = Y 
\end{equation*}
where $Y = (Y_i)_{1 \leq i \leq n}$ with 
\begin{equation}\label{Tan}
Y_i = \lim_{\epsilon \rightarrow 0} \frac{\la s_i,s_i \ra}{2\epsilon}\sum_{w \in W} \int_0^{\cdot} {\bf 1}_{\{0 \leq \la s_i, w^{\star}\Theta_s \ra \leq \epsilon, \Theta_s \in w\overline{C}\}} ds, \quad 1 \leq i \leq n.
\end{equation}
\end{teo}

\begin{note}
1/ Following the terminology used in literature, we refer to $L^0(\Theta)$ as the pushing process of $\pi(\Theta)$. Besides, we shall use the notation $L^0$ for both the pushing process and local times of real semimartingales, nevertheless the argument between brackets indicates which process we refer to. \\
2/ We point the readers interested in Chybiryakov's dissertation to the fact that the author assumed that all roots have equal length $\sqrt{2}$. That is why a rescaling factor showed up in each coordinate $Y_i$.\\
3/ From results derived in \cite{Dai}, one easily sees that the RBM in Weyl chambers is a SRBM in the sense of Williams and Dai, and uniqueness in law holds for the SDE \eqref{SDE}.\\
4/ When $S$ is a basis of $V$, the root system $R$ is said to have full rank. However, this assumption is not always valid as it is not for the root system of type $A$ defined by 
(\cite{Hum} p.41)
\begin{equation*}
R = \{\pm(e_i-e_j), \, 1 \leq i < j \leq n\}, \quad S = \{e_i - e_{i+1}. 1 \leq i \leq n-1\}.
\end{equation*}
As the reader may easily check, the simple system $S$ is no more a basis of $V = \mathbb{R}^n$, rather spans the hyperplane orthogonal to the vector $(1,\dots,1)$ so that $\Lambda$ is a $n \times n-1$ matrix and does not have a left inverse. Nevertheless, the same lines of Chybiryakov's proof lead to 
\begin{equation*}
\Lambda \pi(\Theta_t) = \Lambda \pi(\Theta_0) + \Lambda B_t + Y,
\end{equation*}  
where $Y = (Y_i)_{1 \leq i \leq n-1}$ is defined in the same way as in \eqref{Tan}.  
\end{note}

{\it Review of the Proof}: O. Chybiryakov used a smooth approximation of the projection $\pi$ together with It\^o's formula to derive \eqref{SDE}. But there is an errata in the expression of $\pi$ displayed in \cite{Chy1} p.170. In fact, given $x \in V$, there are as many elements $w$ as the cardinality of the isotropy subgroup of $x$ such that $w^{\star}x \in \overline{C}$. Thus, $w$ is unique if $x \in wC$ (since the isotropy group is trivial, \cite{Hum} p.22) and $\pi = w^{\star}$ and more generally one writes
\begin{align*}
\pi(x) &= \frac{1}{\sharp\{w \in W, w x = x\}} \sum_{w \in W}w^{\star}{\bf 1}_{\{x \in w\overline{C}\}} \\ 
& = \sum_{w \in W}w^{\star}{\bf 1}_{\{x \in wC\}} +  \frac{1}{\sharp\{w \in W, w x = x\}}\sum_{w \in W}w^{\star}{\bf 1}_{\{x \in w(\partial C)\}}
\end{align*}
which is no more differentiable (since the cardinality of the isotropy group is an integer-valued function). However, when one deals with a multidimensional Brownian motion which may hit one and only one hyperplane (\cite{Frie}, p.255), the isotropy group contains two elements: the identity and the reflection orthogonal to that hyperplane thereby  
\begin{align*}
\pi(\Theta_t)  = \sum_{w \in W}w^{\star}{\bf 1}_{\{\Theta_t \in wC\}} +  \frac{1}{2}\sum_{w \in W}w^{\star}{\bf 1}_{\{\Theta_t \in w(\partial C)\}}
\end{align*}
at any time $t$. In this way, Chybiryakov's proof remains valid since $\{t, w^{\star}\Theta_t  \in \partial C\}$ has zero Lebesgue measure by the occupation density formula (\cite{Chy} p.68), and it follows exactly the lines of the one written in \cite{Chy1} p.172-175 with zero multiplicities.

\section{The pushing process as a sum of local times}
Let $s_i \in S$, then it is natural to ask whether $Y_i$ may be written as a sum of the local times $L^0(\la ws_i, \Theta \ra), w \in W$. This claim was shown to be true when $R=B_2$ (\cite{Chy} p.70) but this case is not illustrative since $B_2$ is the most elementary example of a non orthogonal root system. Here we supply the appropriate necessary and sufficient condition ensuring the validity of this claim, thereby recovering Chybiryakov's result and covering other root systems like even dihedral ones (\cite{Dunkl} p.144).  
\begin{pro}
Let $\epsilon > 0$, then $s_i$ is the only simple root in its orbit $R^i = \{ws_i, w \in W\}$ if and only if 
\begin{equation*}
\cup_{w \in W}\{0 \leq  \la ws_i, x \ra \leq \epsilon, x \in w \overline C\} = \cup_{\alpha \in R^i \cap R_+}\{|\la \alpha,x \ra| \leq \epsilon\}
\end{equation*}
for any $x \in V$.
\end{pro}
{\it Proof}: Assume there exists a simple root $s_j, j \neq i$ lying in $R^i$. Then, $s_j^{\perp} \neq s_i^{\perp}$ defines one of the hyperplanes delimiting $\overline{C}$. Since every chamber $w\overline{C}$ is visited only once as $w$ runs over $W$ (counting intersections), then it is impossible to have $0 \leq \la s_j, x\ra \leq \epsilon, x \in \overline{C}$ if one already has $0 \leq \la s_i, x\ra \leq \epsilon, x \in \overline{C}$ since for, $\overline{C}$ will be visited twice. Hence the former set is strictly included into the latter one. Now, assume $s_i$ is the only root in $R^i$, then the boundary of each chamber $w\overline{C}$ contains exactly one hyperplane $H_{\alpha} = \alpha^{\perp}$ for $\alpha = ws_i$.  Since one visits any chamber $wC$ exactly once as $w$ ranges $W$, then one comes close to all hyperplanes $\alpha^{\perp}$ for $\alpha = ws_i \in R^i$, at a distance at most equal to $\epsilon$. This proves that 
\begin{equation*}
\cup_{w \in W}\{0 \leq  \la ws_i, x \ra \leq \epsilon, x \in w \overline C\} = \cup_{\alpha \in R^i}\{0 \leq \la \alpha,x \ra \leq \epsilon\}.
\end{equation*}
Since $-\alpha \in R^i$, the Proposition is proved. $\hfill \blacksquare$  

\begin{cor}
If $s_i$ is the only simple root in its orbit, then 
\begin{equation*}
Y_i = \la s_i,L^0(\Theta) \ra = \sum_{\alpha \in R^i \cap R_+} L^0(\la \alpha,\Theta \ra).
\end{equation*} 
\end{cor}
{\it Proof}: for any $\alpha \in R^i \cap R_+$ and any $\epsilon > 0$, define $A_{\alpha}(x) := \{|\la \alpha,x \ra| \leq \epsilon\}$. Then, one easily shows using the occupation density's formula (Ch.VI \cite{Rev}) that 
\begin{align*}
\sum_{w \in W} \int_0^{\cdot} {\bf 1}_{\{0 \leq \la s_i, w^{\star}\Theta_s \ra \leq \epsilon, \Theta_s \in w\overline{C}\}} ds
 &= \int_0^{\cdot} {\bf 1}_{\cup_{w \in W} \{0 < \la s_i, w^{\star}\Theta_s \ra \leq \epsilon, \Theta_s \in w C\}} ds
 \\& = \int_0^{\cdot} {\bf 1}_{\cup_{\alpha \in R^i \cap R^+} \{A_{\alpha}(\Theta_s) \setminus \la  \alpha,\Theta_s \ra = 0 \ra\}}ds
\\& =  \int_0^{\cdot} {\bf 1}_{\cup_{\alpha \in R^i \cap R^+} \{A_{\alpha}(\Theta_s)\}}ds.
 \end{align*}

Now, one splits the lastly written indicator function into 
\begin{equation*}
\sum_{\alpha \in R^i \cap R^+} {\bf 1}_{A_{\alpha}(\Theta_s)} 
\end{equation*}  
and all the possible intersections of sets $A_{\alpha}(\Theta_s)$, $s$ being fixed. Since $\la \alpha, \alpha\ra = \la s_i,s_i \ra$ for any $\alpha \in R^i$ and since $\la \alpha, \Theta \ra$ is a  local martingale with brackets process $(\la \alpha,\alpha \ra t)_t$, then 
\begin{equation*}
\lim_{\epsilon \rightarrow 0} \frac{\la s_i,s_i \ra}{2\epsilon}\int_0^{\cdot} {\bf 1}_{A_{\alpha}(\Theta_s)}ds = L^0(\la \alpha,\Theta \ra)
\end{equation*}
by Corollary 1.9 Ch.VI in \cite{Rev}. Now, each intersection of sets $A_{\alpha}(\Theta), \alpha \in R^i \cap R_+$ contains two (at least) non necessarily independent Brownian motions (up to scaling factors). Thus, one constructs a planar Brownian motion whose length is less than $c\epsilon$ for some positive constant $c$. As a matter fact, the indicator function of any intersection of sets $A_{\alpha}(\Theta_s), \alpha \in R^i \cap R_+$ is less then 
\begin{equation*}
{\bf 1}_{\{0 \leq Z_s \leq \epsilon\}}
\end{equation*} 
where $Z$ is a Bessel process of dimension $2$. Since the local time of the latter process is almost surely zero, then all integrals of the indicator function of any intersection of sets 
$A_{\alpha}(\Theta), \alpha \in R^i \cap R_+$ vanish after rescaling by $1/\epsilon$ and letting $\epsilon \rightarrow 0$. $\hfill \blacksquare$
  
\begin{nota}
When two simple roots are conjugated, one can only claim 
\begin{equation*}
\cup_{w \in W}\{0 \leq  \la ws_i, x \ra \leq \epsilon, x \in w \overline C\} = \cup_{w \in W}\{-\epsilon \leq  \la ws_i, x \ra \leq \epsilon, x \in w (\overline C \cup s_i\overline{C})\}.
\end{equation*}
An illustrative example is provided by the odd dihedral root system $I_2(3)$ (\cite{Dunkl} p.144): there are six two-dimensional roots represented as
\begin{equation*}
\{\pm e^{-i\pi/2} e^{im\pi/3}, \, m=1,2,3\}
\end{equation*}
and there is only one orbit. The reflections group $W$ contains three rotations of angles $2m\pi/3, m \in \{1,2,3\}$ and three reflections $\mathbb{C} \ni z \mapsto \overline{z}e^{2im\pi/3}$. When choosing $S = \{s_1,s_2\} = \{e^{-i\pi/6}, e^{i\pi/2}\}$, $\overline{C}$ is a wedge of angle $\pi/3$ in the positive quadrant of the plane and 
\begin{equation*}
\cup_{w \in W}\{0 \leq  \la ws_1, x \ra \leq \epsilon, x \in w \overline C\}
\end{equation*} 
is the hachured part of the following picture
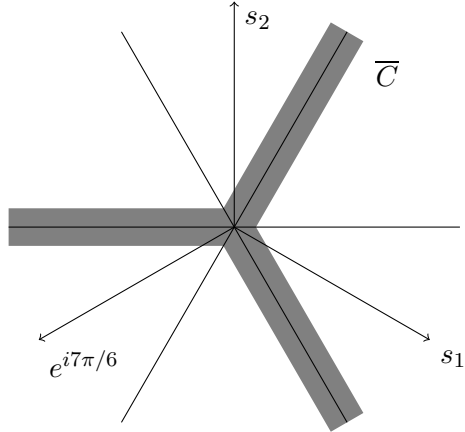
\vspace{0.5cm}
\begin{center}
\begin{tikzpicture}[scale=1]
\draw[line width=5mm, color=gray] (0:0pt) -- (60:3cm) node[at end,below right,color=black]{$\overline{C}$};
\draw[line width=5mm, color=gray] (0:0pt) -- (180:3cm);
\draw[line width=5mm, color=gray] (0:0pt) -- (300:3cm);
\draw[->] (0:0pt) -- (-30:3cm) node[at end,below right]{$s_1$};
\draw[->] (0:0pt) -- (90:3cm) node[at end,below right]{$s_2$};
\draw[->] (0:0pt) -- (-150:3cm) node[at end,below right]{$e^{i7\pi/6}$};
\draw[-] (0:0pt) -- (0:3cm);
\draw[-] (0:0pt) -- (120:3cm);
\draw[-] (0:0pt) -- (240:3cm);
\draw[-] (0:0pt) -- (-60:3cm);
\draw[-] (0:0pt) -- (180:3cm);
\draw[-] (0:0pt) -- (60:3cm);
\end{tikzpicture}
\end{center}
\end{nota}
With some manipulations, one writes $\la L^0(\Theta), s_1\ra$ as
\begin{align*}
\int_0^t {\bf 1}_{\{\la e^{i\pi/3}, \Theta_s \ra \geq 0\}}dL_s^0(\la \Theta, s_1 \ra)&+ \int_0^t {\bf 1}_{\{\la e^{i\pi}, \Theta_s \ra \geq 0\}} dL^0(\la \Theta, s_2\ra)
\\ &+ \int_0^t {\bf 1}_{\{\la e^{i5\pi/3}, \Theta_s \ra \geq 0\}}dL^0(\la \Theta,e^{i7\pi/6} \ra).
\end{align*} 

\section{A multivoque SDE}
Maximal monotone multivoque operators were shown to be a quite efficient tool to deal with SDEs having singular drifts (see references in \cite{Cepa}). They also opened the way to study the strong existence and uniqueness of diffusions that are normally reflected on the boundary of their state spaces (multivoque Skorohod problem). For instance, one deduces from results derived in \cite{Cepa} that there exists a unique strong solution of the following $n$-dimensional SDE: 
\begin{equation}\label{MSDE}
X_t = B_t + n(X_t)dk_t(X) 
\end{equation} 
where $X$ is a continuous process valued in some convex closed domain $\overline{D} \subset \mathbb{R}^n$, $B$ is a $n$-dimensional Brownian motion, $n(x)$ belongs to the inward vector field at $x$ (see below) and $k$ is a continuous process of finite variation such that 
\begin{equation*}
\int_0^t {\bf 1}_{\{X_s \in D\}}d\,|k|_s = 0.
\end{equation*}
In this spirit, a SDE similar to \eqref{MSDE} may be derived for $\pi(\Theta)$ however with oblique reflections. To proceed, recall from \cite{Cepa} that the field of unitary inward normal vectors $n(x)$ to $\overline{C}$ at $x \in \partial C$ is defined by
\begin{equation*}
\{\la n(x), x-a \ra, \, \textrm{for all} \, a \in \overline{C}\}.
\end{equation*}

\begin{pro}
For any hitting time of the boundary $t$, $dL_t^0$ is the sum of an inward normal vector at $\pi(\Theta_t)$ and an inward normal vector at $0_V$. 
\end{pro}

\begin{cor}
The RBM in Weyl chambers satisfies  
\begin{equation*}
\pi(\Theta_t) = \pi(\Theta_0) + B_t + \int_0^t [n(\pi(\Theta_s)) + n_0(s)] dl_s 
\end{equation*}
where $(l_t)_{t \geq 0}$ is an increasing process  and $n_0(s)$ is an inward vector at  $0_V$ at time $s$. 
\end{cor}

{\it Proof}: the expansion of $L^0(\Theta)$ in the dual basis transforms \eqref{SDE} to  
\begin{align*}
\pi(\Theta_t) &= \pi(\Theta_0) + B_t + \sum_{i=1}^n  Y_i  \xi_i. 
\end{align*}
Now expand 
\begin{equation*}
\xi_i  = \sum_{j =1}^n  \la \xi_j,\xi_i \ra s_j
\end{equation*}
and note that any simple root is an inward normal vector at any vactor lying in the corresponding hyperplane, therefore it is so at the vertex $0_V$ since $0_V$ is the intersection of the walls of the simplex. Moreover, $\pi(\Theta)$ cannot visit two or more than two hyperplanes at the same time since $\Theta$ does the same (\cite{Frie}, p.255). Therefore there exists one and only one index $1 \leq i \leq n$ such that $\la \pi(\Theta_t), s_i \ra =0$ at any hitting time of the boundary so that only $Y_i$ survives in the drift at that time. As a matter of fact, it only remains to prove

\begin{lem}\label{Bakry}
$\la \xi_j,\xi_i \ra \geq 0$ for any  $i \neq j$. 
\end{lem}
{\it Proof of the Lemma}: this is a tricky linear algebra exercise and  mainly follows from the fact that $\la s_j,s_i \ra \leq 0$ for any $i \neq j$ (\cite{Hum} p.9). In fact, let $A = (\la s_i,s_j\ra)_{i,j} = \Lambda \Lambda^T$ be the Gram matrix associated with $S$ and let $A^{-1} = (\la \xi_i,\xi_j\ra)_{i,j}$ be its inverse ($\Lambda$ is invertible since $S$ is a basis). Then, since $A$ is symmetric (hence may be diagonalized in an orthogonal basis) and its eigenvalues are strictly positive (it is a positive invertible matrix since $\Lambda$ is invertible), the following functional calculus identity makes perfectly sense: 
\begin{equation*}
A^{-1} = \int_0^{\infty}e^{-tA} dt. 
\end{equation*}  
Since $-A$ has positive offdiagonal entries and negative diagonal ones, then $e^{-tA}$ has positive entries since $e^{-tA} = e^{t(-A + c {\it I})} e^{-ct}$ where $c = \max_{1 \leq i \leq n}
\la s_i,s_i\ra$. The Lemma is then proved. $\hfill \blacksquare$
The Proposition easily follows and the statement of its Corollary is obvious with $l_t =  Y_i(t), n(\pi(\Theta_t)) = \la \xi_i, \xi_i \ra s_i$ and 
\begin{equation*}
n_0(t) = \sum_{j \neq i}\la \xi_i,\xi_j \ra s_j
\end{equation*}
at time $t$ such that $\la \pi(\Theta_t),s_i \ra = 0$.$\hfill \blacksquare$
  
\begin{nota}
The presence of an inward normal vector at $0$ hints to the transience to infinity of $\pi(\Theta)$. Indeed $\Theta$ and $\pi(\Theta)$ has equal length which is a Bessel process of dimension $\geq 2$, the latter process being transient (\cite{Rev}). In particular, when $R = B_2$ (\cite{Hum} p.42), then $\overline{C}$ is a wedge of angle $\pi/4$ and according to the convention adapted in \cite{Var} on the reflections angles (see p.409), these are negative: $\theta_1 = \theta_2 = -\pi/4$ so that the RBM is transient to infinity (see \cite{Will1} p.762).  
\end{nota}

\section{Further developments}
\subsection{Particles system}
By the occupation density's formula (\cite{Rev}), one has 
\begin{equation*}
Y_i = \lim_{\epsilon \rightarrow 0} \frac{\la s_i,s_i \ra}{2\epsilon}\sum_{w \in W} \int_0^{\cdot} {\bf 1}_{\{0 < \la s_i, w^{\star}\Theta_s \ra \leq \epsilon, \Theta_s \in wC\}} ds, \quad 1 \leq i \leq n.
\end{equation*}
Thus, since $\pi(x) = w^{\star}(x)$ when $x \in wC$, then  
\begin{equation*}
Y_i = \frac{\la s_i,s_i \ra}{2}\lim_{\epsilon \rightarrow 0}  \frac{1}{\epsilon}\int_0^{\cdot} {\bf 1}_{\{0 < \la s_i, \pi(\Theta_s) \ra \leq \epsilon\}} ds
\end{equation*}
and again the occupation density formula yields
\begin{equation*}
Y_i = \frac{\la s_i,s_i \ra}{2}\lim_{\epsilon \rightarrow 0}  \frac{1}{\epsilon}\int_0^{\cdot} {\bf 1}_{\{0 \leq \la s_i, \pi(\Theta_s) \ra \leq \epsilon\}} ds.
\end{equation*}
Since $\la s_i, \pi(\Theta) \ra$ is a continuous semimartingale with brakets process $(\la s_i, s_i \ra t)_t$, then (Corollary 1.9, Ch.VI in \cite{Rev}) 
\begin{equation*}
L_t^0(\la s_i, \pi(\Theta) \ra) = 2Y_i
\end{equation*}
where  $L^0(\la s_i, \pi(\Theta) \ra$ is the local time at $0$ of the real-valued semimartingale $\la s_i, \pi(\Theta) \ra$. Hence
\begin{align*}
\pi(\Theta_t) &= \pi(\Theta_0) + B_t + \sum_{i=1}^n  Y_i  \xi_i 
\\& = \pi(\Theta_0) + B_t + \frac{1}{2}\sum_{i=1}^n  L_t^0(\la s_i, \pi(\Theta) \ra) \xi_i. 
\end{align*}

The lastly displayed SDE resembles
\begin{equation*}
dX_t  = dB_t + \sum_{i=1}^n L_t^0(\la n_i, X\ra) V_i 
\end{equation*}
considered in \cite{Sznit} while studying a system of Brownian particles that reflect on hyperplanes orthogonal to $n$-dimensional unit vectors $(n_i)_i$. Above, $V_i$ lies in the hyperplane $n_i^{\perp}$ for all $1 \leq  i \leq n$ so that the framework presented in \cite{Sznit} is slightly different from ours since the RBM in Weyl chambers performs its (oblique) reflections along the dual basis vectors.
 
\subsection{Multidimensional skew Brownian motion} The skew Brownian motion may be defined as a reflected Brownian-like process whose excursions from $0$ change sign according to independent Bernoulli variables of the same parameter $0 \leq p \leq 1$ (see \cite{Lejay} p.434). It is the unique strong solution of (\cite{Lejay} p.431)\footnote{We use a different normalization for the local time.}
\begin{equation}\label{Skew}
X_t  = X_0 + B_t + \frac{(2p-1)}{2}L_t^0(X) 
\end{equation}  
so that one recovers the RBM for $p=1$. Quite complicated ways were proposed to define a multidimensional version of the skew Brownian motion (\cite{Lejay} p.455). Here, one may construct  a continuous process valued in $V$ based on the excursions of $\pi(\Theta)$ from the hyperplanes delimiting $\overline{C}$. Easily speaking, the skew process has continuous paths, evolves like $\pi(\Theta)$ in the interior of chambers $w\overline{C}, w \in W$ and either reflects or transmits at hyperplanes. A stochastic analysis of this process will be the topic of a subsequent paper.

{\bf Acknowledgments}: the author wants to thank Professor D. Bakry for his assistance in proving Lemma \ref{Bakry}. He is also grateful to his colleague J. Zender for providing him with the Latex code used to draw the picture.

\end{document}